\documentclass[12pt]{scrartcl}

\usepackage[T1]{fontenc}
\usepackage[latin1]{inputenc}
\usepackage{amsmath}
\usepackage{amssymb}
\usepackage{amsthm}
\usepackage{dsfont}

\theoremstyle{plain}\newtheorem{theo}{Theorem}
\theoremstyle{plain}\newtheorem{cor}{Corollary}
\theoremstyle{definition}
\theoremstyle{plain}\newtheorem{defi}{Definition}[section]
\theoremstyle{plain}\newtheorem{lem}[defi]{Lemma}
\theoremstyle{plain}\newtheorem{prop}[defi]{Proposition}
\theoremstyle{definition}\newtheorem{exa}[defi]{Example}
\theoremstyle{plain}\newtheorem{assu}{Assumption}

\newcommand{\N}{{\mathds N}}
\newcommand{\R}{{\mathds R}}

\DeclareMathOperator{\cov}{Cov}

\newenvironment{keyword}%
   {\begin{trivlist}\item[]{\bfseries\sffamily Keywords:}\ }
   {\end{trivlist}}
   
\begin{document}

\title{$U$-Processes, $U$-Quantile Processes and Generalized Linear Statistics of Dependent Data}
\author{Martin Wendler\thanks{Martin.Wendler@rub.de} \\ \small Fakult\"{a}t f\"{u}r Mathematik, Ruhr-Universit\"{a}t Bochum, 44780 Bochum, Germany}

\maketitle


\begin{abstract} Generalized linear statistics are an unifying class that contains $U$-statistics, $U$-quantiles, $L$-statistics as well as trimmed and winsorized $U$-statistics. For example, many commonly used estimators of scale fall into this class. $GL$-statistics only have been studied under independence; in this paper, we develop an asymptotic theory for $GL$-statistics of sequences which are strongly mixing or $L^1$ near epoch dependent on an absolutely regular process. For this purpose, we prove an almost sure approximation of the empirical $U$-process by a Gaussian process. With the help of a generalized Bahadur representation, it follows that such a strong invariance principle also holds for the empirical $U$-quantile process and consequently for $GL$-statistics. We obtain central limit theorems and laws of the iterated logarithm for $U$-processes, $U$-quantile processes and $GL$-statistics as straightforward corollaries.

\end{abstract}

\section{Introduction}

\subsection*{$U$-Statistics and the Empirical $U$-Process}

In the whole paper, $\left(X_n\right)_{n\in\N}$ shall be a stationary, real valued sequence of random variables. A $U$-statistic $U_n(g)$ can be described as generalized mean, i.e. the mean of the values $g(X_i,X_j)$, $1\leq i< j\leq n$, where $g$ is a bivariate, symmetric and measurable kernel. The following two estimators of scale are $U$-statistics:

\begin{exa}\label{ex1} Consider $g\left(x,y\right)=\frac{1}{2}\left(x-y\right)^{2}$. A short calculation shows that the related U-statistic is the well-known variance estimator
\begin{equation*}
 U_{n}\left(g\right)=\frac{1}{n-1}\sum_{1\leq i\leq n}\left(X_{i}-\bar{X}\right)^{2}.
\end{equation*}
\end{exa}

\begin{exa}\label{ex2} Let $g\left(x,y\right)=\left|x-y\right|.$ Then the corresponding $U$-statistic is
\begin{equation*}
U_n\left(g\right)=\frac{2}{n(n-1)}\sum_{1\leq i<j\leq n}\left|X_i-X_j\right|,
\end{equation*}
known as Gini's mean difference.
\end{exa}

For $U$-statistics of independent random variables, the central limit theorem (CLT) goes back to Hoeffding \cite{hoef} and was extended to absolutely regular sequences by Yoshihara \cite{yosh}, to near epoch dependent sequences on absolutely regular processes by Denker and Keller \cite{denk} and to strongly mixing random variables by Dehling and Wendler \cite{dehl}. The law of the iterated logarithm (LIL) under independence was proved by Serfling \cite{ser2} and was extended to strongly mixing and near epoch dependent sequences by Dehling and Wendler \cite{deh2}.

Not only $U$-statistics with fixed kernel $g$ are of interest, but also the empirical $U$-distribution function $\left(U_n(t)\right)_{t\in\R}$, which is for fixed $t$ a $U$-statistic with kernel $h(x,y,t):=\mathds{1}_{\left\{g(x,y)\leq t\right\}}$. The Grassberger-Procaccia and the Takens estimator of the correlation dimension in a dynamical system are based on the empirical $U$-distribution function, see Borovkova et al. \cite{boro}.

The functional CLT for the empirical $U$-distribution function has been established by Arcones and Gin\'e \cite{arc4} for independent data, by Arcones and Yu for absolutely regular data \cite{arc3}, and by Borovkova et al. \cite{boro} for data, which is near epoch dependent on absolutely regular processes. The functional LIL for the empirical $U$-distribution function has been proved by Arcones \cite{arc2}, Arcones and Gin\'e \cite{arc5} under independence. The Strong invariance principle has been investigated by Dehling et al. \cite{deh3}. We will show a strong invariance principle under dependence. As a corollary, we will obtain the LIL to sequences which are strongly mixing or $L^1$ near epoch dependent on an absolutely regular process and the CLT under conditions which are slightly different from the conditions in Borovkova et al. \cite{boro}. Let us now proceed with precise definitions:

\begin{defi}\label{def1}
We call a measurable function $h:\R\times\R\times\R\rightarrow\R$, which is symmetric in the first two arguments a kernel function. For fixed $t\in\R$, we call
\begin{equation*}
 U_n\left(t\right):=\frac{2}{n(n-1)}\sum_{1\leq i<j\leq n}h\left(X_i,X_j,t\right)
\end{equation*}
the $U$-statistic with kernel $h\left(\cdot,\cdot,t\right)$ and the process $\left(U_n\left(t\right)\right)_{t\in\R}$ the empirical $U$-distribution function. We define the $U$-distribution function as $U\left(t\right):=E\left[h\left(X,Y,t\right)\right]$, where $X$, $Y$ are independent with the same distribution as $X_1$, and the empirical $U$-process as $\left(\sqrt{n}\left(U_n(t)-U(t)\right)\right)_{t\in\R}$.
\end{defi}

The main tool for the investigation of $U$-statistics is the Hoeffding decomposition into a linear and a so-called degenerate part:
\begin{equation*}
U_n\left(t\right)=U\left(t\right)+\frac{2}{n}\sum_{1\leq i\leq n}h_{1}\left(X_{i},t\right)+\frac{2}{n\left(n-1\right)}\sum_{1\leq i<j\leq n}h_{2}\left(X_{i},X_{j},t\right)
\end{equation*}
where
\begin{align*}
h_1(x,t)&:=Eh(x,Y,t)-U\left(t\right) \\
h_2(x,y,t)&:=h(x,y,t) - h_1(x,t) -h_1(y,t) -U\left(t\right).
\end{align*}

We need some technical assumptions to guarantee the convergence of the empirical $U$-process:

\begin{assu}\label{ass1} The kernel function $h$ is bounded and non-decreasing in the third argument. The $U$-distribution function $U$ is continuous. For all $x,y\in\R$: $\lim_{t\rightarrow\infty}h(x,y,t)=1$, $\lim_{t\rightarrow -\infty}h(x,y,t)=0$.
\end{assu}

Furthermore, we will consider dependent random variables, so we need an additional continuity property of the kernel function (which was introduced by Denker and Keller \cite{denk}):

\begin{assu}\label{ass2} $h$ satisfies the uniform variation condition, that means there is a constant $L$, such that for all $t\in\R$, $\epsilon>0$
\begin{equation*}
 E\left[\sup_{\left\|(x,y)-(X,Y)\right\|\leq \epsilon}\left|h\left(x,y,t\right)-h\left(X,Y,t\right)\right|\right]\leq L\epsilon,
\end{equation*}
where $X$, $Y$ are independent with the same distribution as $X_1$ and $\left\|\cdot\right\|$ denotes the Euclidean norm.
\end{assu}

\subsection*{Empirical $U$-Quantiles and $GL$-Statistics}

For $p\in(0,1)$, the $p$-th $U$-quantile $t_p=U^{-1}(p)$ is the inverse of the $U$-distribution function $U$ at point $p$ (in general, $U$ does not have to be invertible, but this is guaranteed by our Assumption \ref{ass3} at least in the interval $I$ introduced in Theorem \ref{theo2}). A natural estimator of a $U$-quantile is the empirical $U$-quantile $U_n^{-1}(p)$, which is the generalized inverse of the empirical $U$-distribution function at point $p$:
\begin{defi}\label{defi2} Let $p\in(0,1)$ and let $U_n$ be the empirical $U$-distribution function.
\begin{equation*}
U_n^{-1}(p):=\inf\left\{t\big|U_n(t)\geq p\right\}
\end{equation*}
is called the empirical $U$-quantile.
\end{defi}
Empirical $U$-quantiles have applications in robust statistics.

\begin{exa}\label{ex3} Let $h(x,y,t):=\mathds{1}_{\left\{|x-y|\leq t\right\}}$. Then the 0.25-$U$-quantile is the $Q_n$ estimator of scale proposed by Rousseeuw and Croux \cite{rous}, which is highly robust, as its breakdown point is 50\%.
\end{exa}

The kernel function $h(x,y,t):=\mathds{1}_{\left\{|x-y|\leq t\right\}}$ satisfies Assumption \ref{ass2} (uniform variation condition), if the $U$-distribution function is Lipschitz continuous. For every $\epsilon>0$
\begin{multline*}
 E\left[\sup_{\left\|(x,y)-(X,Y)\right\|\leq \epsilon}\left|\mathds{1}_{\left\{|x-y|\leq t\right\}}-\mathds{1}_{\left\{|X-Y|\leq t\right\}}\right|\right]\\
\leq P\left[t-\sqrt{2}\epsilon<|X-Y|\leq t+\sqrt{2}\epsilon\right]\leq U(t+\sqrt{2}\epsilon)-U(t-\sqrt{2}\epsilon) \leq C\epsilon.
\end{multline*}

The empirical $U$-quantile and the empirical $U$-distribution function have a converse behaviour: $U^{-1}_n\left(p\right)$ is greater than $t_p$ iff $U_n\left(t_p\right)$ is smaller than $p$. This motivates a generalized Bahadur representation \cite{baha}:
\begin{equation*}
U^{-1}_n\left(p\right)=t_p+\frac{p-U_n\left(t_p\right)}{u\left(t_p\right)}+R_n(p),
\end{equation*}
where $u=U'$ is the derivative of the $U$-distribution function. For independent data and fixed $p$, Geertsema \cite{geer} established a generalized Bahadur representation with $R_n(p)=O\left(n^{-\frac{3}{4}}\log n\right)$ a.s.. Dehling et al. \cite{deh3} and Choudhury and Serfling \cite{chou} improved the rate to $R_n(p)=O\left(n^{-\frac{3}{4}}(\log n)^{\frac{3}{4}}\right)$. Arcones \cite{arco} proved the exact order $R_n(p)=O\left(n^{-\frac{3}{4}}(\log\log n)^{\frac{3}{4}}\right)$ as for sample quantiles. Under strong mixing and near epoch dependence on an absolutely regular processes, we recently established rates of convergence for $R_n(p)$ which depend on the decrease of the mixing coefficients \cite{wend}. The CLT and the LIL for $U^{-1}_n\left(p\right)$ are straightforward corollaries of the convergence of $R_n$ and the corresponding theorems for $U_n(t_p)$.

In this paper, we will study not a single $U$-quantile, but the empirical $U$-quantile process $\left(U^{-1}_n(p)\right)_{p\in I}$ under dependence, where the interval $I$ is given by $I=[\tilde{C}_1,\tilde{C}_2]$ with $U(C_1)<\tilde{C}_1<\tilde{C}_2<U(C_2)$ and the constants $C_1$, $C_2$ from Assumption \ref{ass3} below. In order to do this, we will examine the rate of convergence of $\sup_{p\in I}R_n(p)$ and use the approximation of the empirical $U$-process by a Gaussian process. As we divide by $u$ in the Bahadur representation, we have to assume that this derivative behaves nicely. Furthermore, we need $U$ to be a bit more than differentiable (but twice differentiable is not needed).

\begin{assu}\label{ass3} $U$ differentiable on an interval $[C_1,C_2]$ with $0<\inf_{t\in[C_1,C_2]} u(t)\leq \sup_{t\in[C_1,C_2]} u(t)<\infty$ ($u(t)=U'(t)$) and
\begin{equation*}
\sup_{t,t'\in[C_1,C_2]:\ \ \left|t-t'\right|\leq x}\left|U(t)-U(t')-u(t)(t-t')\right|=O\left(x^{\frac{5}{4}}\right).
\end{equation*}
\end{assu}

The Bahadur representation for sample quantile process goes back to Kiefer \cite{kief} under independence, Babu and Singh \cite{babu} proved such a representation for mixing data and Kulik \cite{kuli} and Wu \cite{wu} for linear processes, but there seem to be no such results for the $U$-quantile process.

Furthermore, we are interested in linear functionals of the $U$-quantile process.
\begin{defi}\label{defi3} Let $p_1,\ldots,p_d\in I$, $b_1,\dots,b_d\in\R$ and let $J$ be a bounded function, that is continuous a.e. and vanishes outside of $I$. We call a statistic of the form
\begin{multline*} T_n=T\left(U_n^{-1}\right):=\int_{I}J\left(p\right)U_n^{-1}(p)dp+\sum_{j=1}^{d}b_{j}U_n^{-1}(p_{j})\\
=\sum_{i=1}^{\frac{n\left(n-1\right)}{2}}\int_{\frac{2(i-1)}{n\left(n-1\right)}}^{\frac{2i}{n\left(n-1\right)}}J\left(t\right)dt\cdot U_n^{-1}\left(\frac{2i}{n\left(n-1\right)}\right)+\sum_{j=1}^{d}b_{j}U_n^{-1}(p_{j})
\end{multline*}
generalized linear statistic ($GL$-statistic).
\end{defi}
This generalization of $L$-statistics was introduced by Serfling \cite{serf}. $U$-statistics, $U$-quantiles and $L$-statistics can be written as $GL$-statistics (though this might be somewhat artificially). For a $U$-statistics, just take $h(x,y,t)=\mathds{1}_{\left\{g(x,y)\leq t\right\}}$ and $J=1$ (this only works if we can consider the interval $I=[0,1]$). The following example shows how to deal with an ordinary $L$-statistic.

\begin{exa}\label{ex4}Let $h(x,y,t):=\frac{1}{2}\left(\mathds{1}_{\left\{x\leq t\right\}}+\mathds{1}_{\left\{y\leq t\right\}}\right)$, $p_1=0.25$, $p_2=0.75$, $b_1=-1$, $b_2=1$, and $J=0$. Then a short calculation shows that the related $GL$-statistic is  
\begin{equation*}
 T_n=F_n^{-1}(0.75)-F_n^{-1}(0.25),
\end{equation*}
where $F_n^{-1}$ denotes the empirical sample quantile function. This is the well-known inter quartile distance, a robust estimator of scale with 25\% breakdown point. 
\end{exa}

\begin{exa}\label{ex5}Let $h(x,y,t):=\mathds{1}_{\left\{\frac{1}{2}(x-y)^2\leq t\right\}}$, $p_1=0.75$, $b_1=0.25$ and $J(x)=\mathds{1}_{\left\{x\in[0,0.75]\right\}}$. The related $GL$-statistic is called winsorized variance, a robust estimator of scale with 13\% breakdown point. 
\end{exa}
The uniform variation condition also holds in this case, as $h(x,y,t)=\mathds{1}_{\left\{\frac{1}{2}(x-y)^2\leq t\right\}}=\mathds{1}_{\left\{|x-y| \leq \sqrt{2t}\right\}}$ and this is the kernel function of Example \ref{ex3}.

\subsection*{Dependent Sequences of Random Variables}

While the theory of $GL$-statistics under independence has been studied by Serfling \cite{serf}, there seems to be no results under dependence. But many dependent random sequences are very common in applications. Strong mixing and near epoch dependence are widely used concepts to describe short range dependence.
\begin{defi}\label{defi4} Let $\left(X_n\right)_{n\in\N}$ be a stationary process. Then the strong mixing coefficient is given by
\begin{equation*}
\alpha (k) = \sup \left\{\left| P (A\cap B) - P (A) P (B) \right| : A \in \mathcal{F}^n_1, B \in \mathcal{F}^\infty_{n+k}, n \in \N \right\},
\end{equation*}
where $\mathcal{F}^l_a$ is the $\sigma$-field generated by random variables $X_a, \ldots, X_l.$, and $\left(X_n\right)_{n\in\N}$ is called strongly mixing, if $\alpha (k)\rightarrow0$ as $k\rightarrow\infty.$
\end{defi}
Strong mixing in the sense of $\alpha$-mixing is the weakest of the well-known strong mixing conditions, see Bradley \cite{brad}. But this class of weak dependent processes is too strong for many applications, as it excludes examples like linear processes with innovations that do not have a density or data from dynamical systems, see Andrews \cite{andr}.

We will consider sequences which are near epoch dependent on absolutely regular processes, as this class covers linear processes and data from dynamical systems, which are deterministic except for the initial value. Let $T:\left[0,1\right]\rightarrow\left[0,1\right]$ be a piecewise smooth and expanding map such that $\inf_{x\in\left[0,1\right]}\left|T'\left(x\right)\right|>1$. Then there is a stationary process $\left(X_n\right)_{n\in\N}$ such that $X_{n+1}=T\left(X_n\right)$ which can be represented as a functional of an absolutely regular process, for details see Hofbauer and Keller \cite{hofb}. Linear processes (even with discrete innovations) and GARCH processes are also near epoch dependent, see Hansen \cite{hans}. Near epoch dependent random variables are also called approximating functionals (for example in Borovkova et al. \cite{boro}) 

\begin{defi}\label{defi5} Let $\left(X_n\right)_{n\in\N}$ be a stationary process.
\begin{enumerate}
 \item The absolute regularity coefficient is given by
\begin{equation*}
\beta (k) = \sup_{n\in\N}E \sup \{ \left| P (A | \mathcal{F}_{-\infty}^n) - P (A) \right| : A \in \mathcal {F}^\infty_{n + k}\},
\end{equation*}
and $\left(X_n\right)_{n\in\N}$ is called absolutely regular, if $\beta(k)\rightarrow0$ as $k\rightarrow\infty.$
\item We say that $\left(X_n\right)_{n\in\mathds{N}}$ is $L^1$ near epoch dependent on a process $(Z_n)_{n\in\mathds{Z}}$ with approximation constants $(a_l)_{l\in\mathds{N}}$, if
\begin{equation*}\
E \left| X_1 - E (X_1 | \mathcal {G}^l_{- l}) \right| \leq a_l \qquad l = 0, 1,2 \ldots
\end{equation*}
where $\lim_{l \rightarrow \infty} a_l = 0$ and $ \mathcal {G}_{-l}^l$ is the $\sigma$-field generated by $Z_{-l}, \ldots, Z_l.$
\end{enumerate}
\end{defi}

In the literature one often finds $L^2$ near epoch dependence (where the $L^1$ norm in the second part of definition \ref{defi5} is replaced by the $L^2$ norm), but this requires second moments and we are interested in robust estimation. So we want to allow heavier tails and consider $L^1$ near epoch dependence. Furthermore, we do not require that the underlying process is independent, it only has to be weakly dependent in the sense of absolute regularity.

\begin{assu}\label{ass4} Let one of the following two conditions hold:
\begin{enumerate}
 \item $\left(X_n\right)_{n\in\N}$ is strongly mixing with mixing coefficients $\alpha(n)=O(n^{-\alpha})$ for $\alpha\geq 8$ and $E|X_i|^r<\infty$ for a $r>\frac{1}{5}$.
\item $\left(X_n\right)_{n\in\N}$ is near epoch dependent on an absolutely regular process with mixing coefficients $\beta(n)=O(n^{-\beta})$ for $\beta\geq 8$ with approximation constants $a(n)=O(n^{-a})$ for $a=\max\left\{\beta+3,12\right\}$.
\end{enumerate}
\end{assu}

\subsection*{Kiefer-M\"uller processes}

For uniformly on $[0,1]$ distributed and independent random variables $(X_n)_{n\in\N}$, M\"uller \cite{mull} determined the limit distribution of the empirical process \begin{equation*}
\left(\frac{1}{\sqrt{n}}\sum_{1\leq i\leq sn}(\mathds{1}_{\left\{X_i\leq t\right\}}-t)\right)_{t,s\in[0,1]}.
\end{equation*}
It converges weakly towards a Gaussian process $(K(t,s))_{s,t\in[0,1]}$ with covariance function $EK(t,s)K(t',s')=\min\{s,s'\}(\min\{t,t'\}-tt')$. Kiefer \cite{kie2} proved an almost sure invariance principle: After enlarging the probability space, there exists a copy of the Kiefer-M\"uller process $K$ such that the empirical process and $K$ are close together with respect to the supremum norm. Berkes and Philipp \cite{berk} extended this to dependent random variables. For sample quantiles, Cs\"org\H{o} and R\'ev\'esz \cite{csor} established a strong invariance principle, but only under independence. We will extend this to dependent data and to $U$-quantiles.

A strong invariance principle is a very interesting asymptotic theorem, as the limit behaviour of Gaussian processes is well understood and it is then possible to conclude that the approximated process has the same asymptotic properties. Note that a Kiefer-M\"uller processes can be described as a functional Brownian motion, as its increments in $s$ direction are independent Brownian Bridges. We have the following scaling behaviour: $(\frac{1}{\sqrt{n}}K(t,ns))_{s,t\in[0,1]}$ has the same distribution as $(K(t,s))_{s,t\in[0,1]}$.

Furthermore, a functional LIL holds: The sequence
\begin{equation*}
\left((\frac{1}{\sqrt{2n\log\log n}}K(t,ns))_{s,t\in[0,1]}\right)_{n\in\N}
\end{equation*}
is almost surely relatively compact (with respect to the supremum norm). The limit set is the unit ball of the reproducing kernel Hilbert space associated with the covariance function of the process $(K(t,s))_{s,t\in[0,1]}$. For details about the reproducing kernel Hilbert space, see Aronszajn \cite{aron} or Lai \cite{lai}.

\section{Main Results}

\subsection*{Empirical $U$-Process}

The asymptotic theory for the empirical $U$-process makes use of the Hoeffding decomposition, recall that $h_1(x,t):=E\left[h(x,Y,t)\right]-U(t)$. Under Assumptions \ref{ass1}, \ref{ass2} and \ref{ass4}, the following covariance function converges absolutely and is continuous (compare to Theorem 5 of Borovkova et al. \cite{boro}):
\begin{multline*}
\Gamma(t,t')=4\cov\left[h_1\left(X_1,t\right),h_1\left(X_1,t'\right)\right]\\
+4\sum_{k=1}^{\infty}\cov\left[h_1\left(X_1,t\right),h_1\left(X_{k+1},t'\right)\right]+4\sum_{k=1}^{\infty}\cov\left[h_1\left(X_{k+1},t\right),h_1\left(X_1,t'\right)\right].
\end{multline*}

\begin{theo}\label{theo1}
Under the assumptions \ref{ass1}, \ref{ass2} and \ref{ass4} there exists a centered Gaussian process $(K(t,s))_{t,s\in\R}$ (after enlarging the probability space if necessary) with covariance function
\begin{equation*}
EK(t,s)K(t',s')=\min\left\{s,s'\right\}\Gamma(t,t')
\end{equation*}
such that almost surely
\begin{equation*}
\sup_{\substack{t\in\R\\s\in[0,1]}}\frac{1}{\sqrt{n}}\left|\lfloor ns\rfloor (U_{\lfloor ns\rfloor}(t)-U(t))-K(t,ns)\right|=O(\log^{-\frac{1}{3840}}n).
\end{equation*}
\end{theo}

The rate of convergence to zero in this theorem is very slow, but the same as in Berkes and Philipp \cite{berk}, as we strongly use their method of proof. By the scaling property of the process $K$, we obtain the asymptotic distribution of $U_{\lfloor ns\rfloor}(t)$, and by Theorem 2.3 of Arcones \cite{arc6} a functional LIL:

\begin{cor}\label{cor1}
Under the assumptions \ref{ass1}, \ref{ass2} and \ref{ass4} the empirical $U$-process
\begin{equation*}
 \left(\frac{\lfloor ns\rfloor }{\sqrt{n}}(U_{\lfloor ns\rfloor}(t)-U(t))\right)_{t\in\R,s\in[0,1]}
\end{equation*}
converges weakly in the space $D(\R\times[0,1])$ (equipped with the supremum norm) to a centered Gaussian Process $(K(t,s))_{t,s\in\R}$ introduced in Theorem \ref{theo1}. The sequence
\begin{equation*}
 \left(\left(\frac{\lfloor ns\rfloor }{\sqrt{2n\log\log n}}(U_{\lfloor ns\rfloor}(t)-U(t))\right)_{t\in\R,s\in[0,1]}\right)_{n\in\N}
\end{equation*}
is almost surely relatively compact in the space $D(\R\times[0,1])$ (equipped with the supremum norm) and the limit set is the unit ball of the reproducing kernel Hilbert space associated with the covariance function of the process $K$.
\end{cor}

The first part of this corollary is very similar to Theorem 9 of Borovkova et al. \cite{boro} (they use a continuity condition that is different from our Assumption \ref{ass2}). Up to our knowledge, part 2 is the first functional LIL for empirical $U$-processes under dependence.

\subsection*{Generalized Bahadur Representation}
Recall that the remainder term in the generalized Bahadur representation is defined as
\begin{equation*}
 R_n(p)=U_n^{-1}\left(p\right)-t_p-\frac{p-U_n\left(t_p\right)}{u\left(t_p\right)}
\end{equation*}
and that we write $t_p:=U^{-1}(p)$. We set $U_0^{-1}(p):=0$ as it is not possible to find a generalized inverse of $U_0=0$.
\begin{theo}\label{theo2} Under the Assumptions \ref{ass1}, \ref{ass2}, \ref{ass3} and \ref{ass4}
\begin{equation*}
\sup_{\substack{p\in I\\s\in[0,1]}}\frac{\lfloor ns\rfloor}{\sqrt{n}}|R_{\lfloor ns\rfloor}(p)|=o(n^{-\frac{\gamma}{8}}\log n)
\end{equation*}
almost surely with $I=[\tilde{C}_1,\tilde{C}_2]$, where $U(C_1)<\tilde{C}_1<\tilde{C}_2<U(C_2)$, $\gamma:=\frac{\alpha-2}{\alpha}$ (if the first part of Assumption \ref{ass4} holds) respectively $\gamma:=\frac{\beta-3}{\beta+1}$ (if the second part of Assumption \ref{ass4} holds).
\end{theo}

Note that for a fast decay of the mixing coefficients, the rate becomes close to $n^{-\frac{1}{8}}$, while the optimal rate for sample quantiles of independent data is $n^{-\frac{1}{4}}(\log n)^{\frac{1}{2}}(\log \log n)^{\frac{1}{4}}$.

\subsection*{Empirical $U$-Quantiles and $GL$-Statistics}

Using the Bahadur representation, we can deduce the asymptotic behaviour of the empirical $U$-quantile process from Theorem \ref{theo1}.

\begin{theo}\label{theo3} Under the Assumptions \ref{ass1}, \ref{ass2}, \ref{ass3} and \ref{ass4}, there exists a centered Gaussian process $(K'(p,s))_{p\in I,s\in\R}$ (after enlarging the probability space if necessary), where $I$ is the interval introduced in Theorem \ref{theo2}, with covariance function
\begin{equation*}
EK'(p,s)K'(p',s')=\min\left\{s,s'\right\}\frac{1}{u(t_p)u(t_{p'})}\Gamma(t_p,t_{p'})
\end{equation*}
such that
\begin{equation*}
\sup_{\substack{p\in I\\s\in[0,1]}}\frac{1}{\sqrt{n}}\left|\lfloor ns\rfloor (U^{-1}_{\lfloor ns\rfloor}(p)-t_p)-K'(p,ns)\right|=O(\log^{-\frac{1}{3840}}n).
\end{equation*}
\end{theo}

$K'$ is a Gaussian process with independent increments in $s$ direction, so we have the following consequences:

\begin{cor}\label{cor2} Under the Assumptions \ref{ass1}, \ref{ass2}, \ref{ass3} and \ref{ass4}
\begin{equation*}
 \left(\frac{\lfloor ns\rfloor }{\sqrt{n}}(U^{-1}_{\lfloor ns\rfloor}(p)-t_p)\right)_{p\in I,s\in[0,1]}
\end{equation*}
converges weakly in the space $D(I\times[0,1])$ (equipped with the supremum norm) to the centered Gaussian Process $(K'(p,s))_{p\in I,s\in\R}$ introduced in Theorem \ref{theo3}.
The sequence
\begin{equation*}
 \left(\left(\frac{\lfloor ns\rfloor }{\sqrt{2n\log\log n}}(U^{-1}_{\lfloor ns\rfloor}(p)-t_p)\right)_{p\in I,s\in[0,1]}\right)_{n\in\N}
\end{equation*} 
is almost surely relatively compact in the space $D(I\times[0,1])$ (equipped with the supremum norm) and the limit set is the unit ball of the reproducing kernel Hilbert space associated with the covariance function of the process $K'$.
\end{cor}

As $GL$-statistics are linear functionals of the empirical $U$-quantile process, we get an approximation for $T_n$:

\begin{theo}\label{theo4}Let $p_1,\ldots,p_d\in I$ and let $J$ be a bounded function, that is continuous a.e. and vanishes outside of $I$. Under the assumptions \ref{ass1}, \ref{ass2}, \ref{ass3} and \ref{ass4}, there exists (after enlarging the probability space if necessary) a Brownian motion $B$, such that for $T_n$ defined in Definition \ref{defi3} and
\begin{multline*}
\sigma^2=\int_{\tilde{C_1}}^{\tilde{C_2}}\int_{\tilde{C_1}}^{\tilde{C_2}}\frac{\Gamma(t_p,t_q)}{u(t_p)u(t_q)}J(p)J(q)dpdq\\
+2\sum_{j=1}^db_j\int_{\tilde{C_1}}^{\tilde{C_2}}\frac{\Gamma(t_{p_j},t_p)}{u(t_{p_j})u(t_p)}J(p)dp+\sum_{i,j=1}^db_ib_j\frac{\Gamma(t_{p_i},t_{p_j})}{u(t_{p_i})u(t_{p_j})}
\end{multline*}
we have that
\begin{equation*}
\sup_{s\in[0,1]}\frac{1}{\sqrt{n}}\left|\lfloor ns\rfloor (T_{\lfloor ns\rfloor}-T(U^{-1}))-\sigma B(ns)\right|=O(\log^{-\frac{1}{3840}}n)
\end{equation*}
almost surely.
\end{theo}

By the well-known properties of Brownian motions, we have:

\begin{cor}\label{cor3}Let $p_1,\ldots,p_d\in I$ and let $J$ be a bounded function. Under the assumptions \ref{ass1}, \ref{ass2}, \ref{ass3} and \ref{ass4} for $T_n$ defined in Definition \ref{defi3}:
\begin{equation*}
\frac{\lfloor ns\rfloor}{\sqrt{n}} (T_{\lfloor ns\rfloor}-T(U^{-1}))
\end{equation*}
converges weakly to the Brownian motion $\sigma B(s)$ with $\sigma^2$ as in Theorem \ref{theo4}. Furthermore, we have that the sequence
\begin{equation*}
\left(\frac{\lfloor ns\rfloor}{\sqrt{2n\log\log n}} (T_{\lfloor ns\rfloor}-T(U^{-1}))_{s\in[0,1]}\right)_{n\in\N}
\end{equation*}
is almost surely relatively compact in the space of bounded continuous functions $C[0,1]$ (equipped with the supremum norm) and the limit set is
\begin{equation*}
\left\{f:[0,1]\rightarrow\R\big|f(0)=0,\ \int_0^1 f'^2(s)ds\leq \sigma^2 \right\}.
\end{equation*}
\end{cor}

\section{Preliminary Results}

\begin{prop}\label{pro1}
Under the assumptions \ref{ass1}, \ref{ass2} and \ref{ass4} there exists a centered Gaussian process $(K(t,s))_{t,s\in\R}$ (after enlarging the probability space if necessary) with covariance function
\begin{equation*}
EK(t,s)K(t',s')=\min\left\{s,s'\right\}\Gamma(t,t')
\end{equation*}
such that almost surely
\begin{equation*}
\sup_{\substack{t\in\R\\s\in[0,1]}}\frac{1}{\sqrt{n}}\left|\left(2\sum_{1\leq i\leq ns}h_1(X_i,t)-K(t,ns)\right)\right|=O(\log^{-\frac{1}{3840}}n).
\end{equation*}
\end{prop}

\begin{proof} This proposition is basically Theorem 1 of Berkes and Philipp \cite{berk}, which we have to generalize in three aspects: 
\begin{enumerate}
\item Berkes and Philipp assume that the covariance kernel $\Gamma$ is positive definite, we want to avoid this condition here.
\item Berkes and Philipp consider indicator functions $\mathds{1}_{\left\{x\leq t\right\}}$, while in this version of the proposition, we deal with more general functions $Eh(x,Y,t)$.
\item Theorem 1 of Berkes and Philipp is restricted to the distribution function $F(t)=E\mathds{1}_{\left\{X_i\leq t\right\}}=t$, we will extend this to a function $U$ according to our Assumption \ref{ass1}.
\end{enumerate}
The mixing condition of Berkes and Philipp is the same as our Assumption \ref{ass4}.
\begin{enumerate}
\item In the proof of their Theorem 1, Berkes and Philipp use the fact that $\Gamma$ is positive definite only for two steps. Their Proposition 4.1 (page 124) also holds if this is not the case. It is easy to see that the characteristic functions of the finite dimensional distributions then might converge to 1 at some points, but with the required rate. Furthermore, we have to show (page 135) that for all $t_{1},\ldots,t_{d_k}\in[0,1]$, $P[\|(K(t_1,1),\ldots,K(t_{d_k},1))\|\geq\frac{1}{4}T_k]\leq\delta_k$, where $T_k$ and $\delta_k$ are defined in their article. Let $\Gamma_{d_k}=\left(\Gamma(t_i,t_j)\right)_{1\leq i,j\leq d_k}$ be the covariance matrix of $K(t_1,1),\ldots,K(t_{d_k},1)$ and $\rho$ its biggest eigenvalue. We first consider the case that $\rho>0$. As $\Gamma_{d_k}$ is symmetric and positive semidefinite, there exist a matrix $\Gamma^{\frac{1}{2}}_{d_k}$ such that $\left(\Gamma^{\frac{1}{2}}_{d_k}\right)^t\Gamma^{\frac{1}{2}}_{d_k}=\Gamma_{d_k}$ and the vector $K(t_1,1),\ldots,K(t_{d_k},1)$ has the same distribution as $\Gamma^{\frac{1}{2}}_{d_k}(W_1,\ldots,W_{d_k})^t$, where $W_1,\ldots,W_{d_k}$ are independent standard normal random variables. So it follows that
\begin{multline*}
P[\|(K(t_1,1),\ldots,K(t_{d_k}))\|\geq\frac{1}{4}T_k]=P[\|\Gamma^{\frac{1}{2}}_{d_k}(W_1,\ldots,W_{d_k})\|\geq\frac{1}{4}T_k]\\
\leq P[\sqrt{\rho}\|(W_1,\ldots,W_{d_k})^t\|\geq\frac{1}{4}T_k]\\
=\frac{1}{(2\pi)^{\frac{1}{2}d_{k}}}\int_{\|(x_1,\ldots,x_{d_k})\|\geq\frac{1}{4\sqrt{\rho}}T_k}\exp(-\frac{1}{2}(x_1^2+\ldots+x_{d_k}^2))dx_1\ldots dx_{d_k}.
\end{multline*}
The rest of the proof is then exactly the same as in Berkes and Philipp \cite{berk}. In the case $\rho=0$, we have that $\Gamma=0$, so trivially $P[\|(K(t_1,1),\ldots,K(t_{d_k}))\|\geq\frac{1}{4}T_k]=0\leq \delta_k$.
\item The proof uses different properties of the indicator functions. If the process $(X_n)_{n\in\N}$ is near epoch dependent with constants $(a_n)_{n\in\N}$, then as a consequence of Lemma 3.2.1 of Philipp \cite{phil} the process $\left(\mathds{1}_{\left\{X_n\leq t\right\}}\right)_{n\in\N}$ is near epoch dependent with constants $(\sqrt{a_n})_{n\in\N}$. The same holds for the sequence $(h_1(X_n,t))_{n\in\N}$ by Assumption \ref{ass2}, Lemma 3.5 and 3.10 of Wendler \cite{wend}.

Furthermore, $h$ and $U$ are non-decreasing in $t$. Berkes and Philipp used different moment properties, which we also assume: $h_1(X_n,t)$ is bounded by 1 and $E|h_1(X_n,t)-h_1(X_n,t')|\leq C|t-t'|$ for $t,t'\in\R$, so consequently for $m\geq1$ $\left\|h_1(X_n,t)\right\|_m\leq 1$ and $\left\|h_1(X_n,t)-h_1(X_n,t')\right\|_m\leq|t-t'|^{\frac{1}{m}}$. So this more general version can be proved along the lines of the proof in Berkes and Philipp \cite{berk}.
\item If $U(t)=t$ does not hold, note that $Eh_1(X_i,t_p)=U(t_p)=p$ with $t_p=U^{-1}(p):=\inf\{t\in\R|U(t)\geq p\}$, because $U$ is continuous. Clearly, Assumption \ref{ass1} and \ref{ass2} hold for $h(x,y,U^{-1}(p))$. Furthermore, notice that if $U(t)=U(s)$, we have that $h_1(X_i,t)=h_1(X_i,s)$ almost surely by monotonicity of $h$, so
\begin{equation*}
\sum_{i=1}^{n}h_1(X_i,t)=\sum_{i=1}^{n}h_1(X_i,t_{U(t)})
\end{equation*}
almost surely. From the first two parts of the proof, we know that there is a centered Gaussian process $K^\star$ with covariance function
\begin{equation*}
E[K^\star(p,s)K^\star(p',s')]=\min\left\{s,s'\right\}\Gamma(t_p,t_{p'})
\end{equation*}
with
\begin{equation*}
\sup_{\substack{p\in[0,1]\\s\in[0,1]}}\frac{1}{\sqrt{n}}\left|\left(2\sum_{1\leq i\leq ns}h_1(X_i,t_p)-K^\star(p,ns)\right)\right|=O(\log^{-\frac{1}{3840}}n).
\end{equation*}
almost surely. The Gaussian process $K$ with $K(t,s)=K^\star(U(t),s)$ has the required covariance function and
\begin{multline*}
\sup_{\substack{t\in\R\\s\in[0,1]}}\frac{1}{\sqrt{n}}\left|\left(2\sum_{1\leq i\leq ns}h_1(X_i,t)-K(t,ns)\right)\right|=\\
\sup_{\substack{t\in\R\\s\in[0,1]}}\frac{1}{\sqrt{n}}\left|\left(2\sum_{1\leq i\leq ns}h_1(X_i,t_{U(t)})-K^\star(U(t),ns)\right)\right|=O(\log^{-\frac{1}{3840}}n).
\end{multline*}
\end{enumerate}
\end{proof}

\begin{lem}\label{lem1} Let $C_3,C_4,L$ be positive constants. Under Assumption \ref{ass4}, part 1 (strong mixing), there exists a constant $C$, such that for all measurable, non-negative functions $g:\R\rightarrow\R$ that are bounded by $C_3$ with $E\left|g\left(X_1\right)-Eg\left(X_1\right)\right|\geq C_4n^{-\frac{\alpha}{\alpha+1}}$ and satisfy the variation condition with constant $L$, and all $n\in\N$ we have
\begin{equation*}
 E\left(\sum_{i=1}^{n}g\left(X_i\right)-E\left[g\left(X_1\right)\right]\right)^4\leq C n^2\left(\log n\right)^2\left(E\left|g\left(X_1\right)\right|\right)^{1+\gamma},
\end{equation*}
where $\gamma$ is defined in Theorem \ref{theo2}. The same statement holds under Assumption \ref{ass4}, part 2 (near epoch dependence on absolutely regular sequence) for functions $g:\R\rightarrow\R$ with $E\left|g\left(X_1\right)-Eg\left(X_1\right)\right|\geq C_4n^{-\frac{\beta}{\beta+1}}$.
\end{lem}

This is Lemma 3.4 respectively 3.6 of Wendler \cite{wend}.

\begin{lem}\label{lem2}
Under Assumptions \ref{ass1}, \ref{ass2} and \ref{ass4}, there exists a constant $C$, such that for all $t\in\R$ and all $n\in\N$
\begin{equation*}
 \sum_{i_1,j_1,i_2,j_2=1}^n \left|E\left[h_2(X_{i_1},X_{j_1},t)h_2(X_{i_2},X_{j_2},t)\right]\right|\leq Cn^2.
\end{equation*}
\end{lem}
This is Lemma 4.4 of Dehling and Wendler \cite{deh2}.

\begin{lem}\label{lem3} Under the Assumptions \ref{ass1}, \ref{ass2} and \ref{ass4}
\begin{equation*}
\sup_{t\in\R}\left|\sum_{1\leq i<j\leq n}h_{2}\left(X_{i},X_{j},t\right)\right|=o\left(n^{\frac{3}{2}-\frac{\gamma}{8}}\right)
\end{equation*}
almost surely with $\gamma$ as in Theorem \ref{theo2}.
\end{lem}

In all our proofs, $C$ denotes a constant and may have different values from line to line.

\begin{proof} Without loss of generality, we can assume that $U(t)=t$, otherwise we use the same transformation as in the proof of Proposition \ref{pro1} and study the kernel function $h(x,y,U^{-1}(p))$. We define $Q_n(t):=\sum_{1\leq i<j\leq n}h_{2}\left(X_{i},X_{j},t\right)$. For $l\in\N$, let $k=k_l=2^{\lceil\frac{5}{8}l\rceil}$ and $t_{r,l}=\frac{r}{k_l}$ for $r=0,\ldots,k_l$, so that $Cn^{\frac{5}{8}}\leq t_{r,l}-t_{r-1,l}=\frac{1}{k_l}\leq C'n^{\frac{5}{8}}$ for all $n\in\N$ with $2^{l-1}\leq n<2^{l}$ and some constants $C,C'$. By Assumption \ref{ass1}, $h$ and $U$ are non-decreasing in $t$, so we have for any $t\in[t_{r-1,l},t_{r,l}]$, $ n<2^{l}$
\begin{multline*}
\left|Q_n(t)\right|=\left|\sum_{1\leq i<j\leq n}\left(h\left(X_{i},X_{j},t\right)-h_1(X_{i},t)-h_1(X_{j},t))-U(t)\right)\right|\\
\leq \max\left\{\left|\sum_{1\leq i<j\leq n}\left(h\left(X_{i},X_{j},t_{r,l}\right)-h_1(X_{i},t)-h_1(X_{j},t)-U(t)\right)\right|\right.,\\
\left.\left|\sum_{1\leq i<j\leq n}\left(h\left(X_{i},X_{j},t_{r-1,l}\right)-h_1(X_{i},t)-h_1(X_{j},t)-U(t)\right)\right|\right\}\\
\leq \max\left\{|Q_{n}(t_{r,l})|,|Q_{n}(t_{r-1,l})|\right\}\\
+(n-1)\max\left\{\left|\sum_{i=1}^{n}(h_1(X_{i},t_{r,l})-h_1(X_{i},t)))\right|,\left|\sum_{i=1}^{n}(h_1(X_{i},t)-h_1(X_{i},t_{r-1,l})))\right|\right\}\\
+\frac{n(n-1)}{2}|U(t_{r,l})-U(t_{r-1,l})|\\
\leq \max\left\{|Q_n(t_{r,l})|,|Q_n(t_{r-1,l})|\right\}\\
+(n-1)\left|\sum_{i=1}^n(h_1(X_{i},t_{r,l})-h_1(X_{i},t_{r-1,l})))\right|+2\frac{n(n-1)}{2}|U(t_{r,l})-U(t_{r-1,l})|.
\end{multline*}
So we have that
\begin{multline*}
 \sup_{t\in\R}\left|Q_{n}(t)\right|\\
\leq \max_{r=0,\ldots,k}\left|Q_{n}(t_{r,l})\right|+\max_{r=0,\ldots,k}(n-1)\left|\sum_{i=1}^{n}\left(h_1(X_{i},t_{r,l})-h_1(X_{i},t_{r-1,l})\right)\right|\\
+\max_{r=0,\ldots,k}n(n-1)|U(t_{r,l})-U(t_{r-1,l})|.
\end{multline*}
We will treat these three summands separately. By the choice of $t_1,\ldots,t_{k-1}$, we have $|U(t_{r,l})-U(t_{r-1,l})|=t_{r,l}-t_{r-1,l}=\frac{1}{k_l}$, so for the last summand and $2^{l-1}\leq n<2^{l}$ we know that $\max_{r=0,\ldots,k}n(n-1)|U(t_{r,l})-U(t_{r-1,l})|\leq Cn^{2-\frac{5}{8}}=o\left(n^{\frac{3}{2}-\frac{\gamma}{8}}\right)$ . For the first summand, we obtain by similar arguments as the ones used by Wu \cite{wu2} to prove his inequality (6) of his Proposition 1 or by Dehling and Wendler \cite{deh2} to prove their line (5)
\begin{multline*}
E[\max_{n=2^{l-1},\ldots,2^{l}-1}\max_{r=0,\ldots,k}\left|Q_n(t_{r,l})\right|^2]\\
\leq \sum_{r=0}^kE\left[\left(\sum_{d=1}^{l}\max_{i=1,\ldots,2^{l-d}}\left|Q_{i2^{d-1}}(t_{r,l})-Q_{(i-1)2^{d-1}}(t_{r,l})\right|\right)^2\right]\\
\leq \sum_{r=0}^kl\sum_{d=1}^{l}\sum_{i=1}^{2^{l-d}}E\left[\left(Q_{i2^{d-1}}(t_{r,l})-Q_{(i-1)2^{d-1}}(t_{r,l})\right)^2\right]\\
\leq \sum_{r=0}^kl\sum_{d=1}^{l}\sum_{i_1,j_1,i_2,j_2=1}^{2^{l}} \left|E\left[h_2(X_{i_1},X_{j_1},t)h_2(X_{i_2},X_{j_2},t)\right]\right|\\
\leq Ckl^2 2^{2(l+1)}\leq Cl^2 2^{(2+\frac{5}{8})l},
\end{multline*}
where we used Lemma \ref{lem2} in the last line. With the Chebyshev inequality, it follows for every $\epsilon>0$
\begin{multline*}
 \sum_{l=1}^\infty P\left[\max_{n=2^{l-1},\ldots,2^{l}-1}\max_{r=0,\ldots,k}\left|Q_n(t_{r,l})\right|>\epsilon 2^{l(\frac{3}{2}-\frac{\gamma}{8})}\right]\\
\leq\sum_{l=1}^\infty \frac{1}{\epsilon^22^{l(3-\frac{\gamma}{4})}}E[\max_{n=1,\ldots,2^{l}}\max_{r=0,\ldots,k}\left|Q_n(t_{r,l})\right|^2]\leq \sum_{l=1}^\infty\frac{1}{\epsilon^22^{l(3-\frac{\gamma}{4})}}l^2 2^{(2+\frac{5}{8})l}<\infty,
\end{multline*}
as $\gamma\leq 1$, so by the Borel Cantelli lemma
\begin{equation*}
P\left[\max_{n=2^{l-1},\ldots,2^{l}-1}\max_{r=0,\ldots,k}\left|Q_n(t_{r,l})\right|>\epsilon 2^{l(\frac{3}{2}-\frac{\gamma}{8})}\ \text{i.o.}\right]=0
\end{equation*}
(the meaning of the abbreviation i.o. is ``infinitely often''). It remains to show the convergence of the second summand:
\begin{multline*}
E\left(\max_{n=2^{l-1},\ldots,2^{l}-1}\max_{r=1,\ldots,k}(n-1)\left|\sum_{i=1}^n(h_1(X_{i},t_{r,l})-h_1(X_{i},t_{r-1,l}))\right|\right)^4\\
\leq
2^{4(l+1)}\sum_{r=1}^kE\left(\max_{n=2^{l-1},\ldots,2^{l}-1}\left|\sum_{i=1}^n(h_1(X_{i},t_{r,l})-h_1(X_{i},t_{r-1,l}))\right|\right)^4\\
\leq C2^{6l} l^2k(\max_{r=1,\ldots,k}|t_{r,l}-t_{r-1,l}|)^{1+\gamma}\leq Cl^22^{(6-\frac{5}{8}\gamma)l},
\end{multline*}
where we used Corollary 1 of M\'oricz and Lemma \ref{lem1} to obtain the last line. Remember that $k=k_l=O\left( 2^{\frac{5}{8}l}\right)$ and that $\left|t_{r,l}-t_{r-1,l}\right|\geq\frac{1}{2^{\frac{5}{8}l}}$. We conclude that
\begin{multline*}
\sum_{l=0}^\infty P\left[\max_{n=2^{l-1},\ldots,2^{l}-1}\max_{r=1,\ldots,k}(n-1)\left|\sum_{i=1}^n(h_1(X_{i},t_{r,l})-h_1(X_{i},t_{r-1,l}))\right|>\epsilon 2^{(\frac{3}{2}-\frac{\gamma}{8})l}\right]\\
\leq \sum_{l=0}^\infty\frac{C}{\epsilon^42^{l(6-\frac{\gamma}{2})}}E\left(\max_{n=2^{l-1},\ldots,2^{l}-1}\max_{r=1,\ldots,k}(n-1)\left|\sum_{i=1}^n(h_1(X_{i},t_{r,l})-h_1(X_{i},t_{r-1,l}))\right|\right)^4\\
\leq \sum_{l=0}^\infty\frac{C}{\epsilon^42^{l(6-\frac{\gamma}{2})}}l^22^{(6-\frac{5}{8}\gamma)l}=\sum_{l=0}^\infty\frac{Cl^2}{\epsilon^42^{\frac{\gamma}{8}l}}<\infty.
\end{multline*}
The Borel Cantelli lemma completes the proof.
\end{proof}

\begin{lem}\label{lem4} Let $F$ be a non-decreasing function, $c,l>0$ constants and $[C_1,C_2]\subset\R$. If for all $t,t'\in[C_1,C_2]$ with $|t-t'|\leq l+2c$
\begin{equation*}
|F(t)-F(t')-(t-t')|\leq c,
\end{equation*}
then for all $p,p'\in\R$ with $|p-p'|\leq l$ and $F^{-1}(p),F^{-1}(p')\in(C_1+2c+l,C_2-2c-l)$
\begin{equation*}
|F^{-1}(p)-F^{-1}(p')-(p-p')|\leq c
\end{equation*}
where $F^{-1}(p):=\inf\left\{t\big|F(t)\geq p\right\}$ is the generalized inverse.
\end{lem}

\begin{proof} Without loss of generality we assume that $p<p'$. Let $\epsilon\in(0,c)$. By our assumptions
\begin{multline*}
F\left(F^{-1}(p)+(p'-p)+c+\epsilon\right)\geq F\left(F^{-1}(p)+\epsilon\right)+(p'-p)+c-c\\
\geq p+(p'-p)=p'.
\end{multline*}
By the definition of $F^{-1}$, it follows that
\begin{equation*}
F^{-1}(p')=\inf\left\{t\big|F(t)\geq p'\right\}\leq F^{-1}(p)+(p'-p)+c+\epsilon.
\end{equation*}
So taking the limit $\epsilon\rightarrow0$, we obtain
\begin{equation*}
F^{-1}(p')\leq F^{-1}(p)+(p'-p)+c.
\end{equation*}
On the other hand
\begin{multline*}
F\left(F^{-1}(p)+(p'-p)-c-\epsilon\right)\leq F\left(F^{-1}(p)-\epsilon\right)+(p'-p)-c+c\\
\leq p+(p'-p)=p'.
\end{multline*}
So we have that
\begin{equation*}
F^{-1}(p')\geq F^{-1}(p)+(p'-p)-c-\epsilon,
\end{equation*}
and hence $F^{-1}(p')\geq F^{-1}(p)+(p'-p)-c$. Combining the upper and lower inequality for $F^{-1}(p')$, we conclude that $|F^{-1}(p)-F^{-1}(p')-(p-p')|\leq c$.
\end{proof}

\begin{lem}\label{lem5} Under the Assumptions \ref{ass1}, \ref{ass2}, \ref{ass3} and \ref{ass4} for any constand $C>0$
\begin{equation*}
 \sup_{\substack{t,t'\in[C_1,C_2]:\\ |t-t'|\leq C\sqrt{\frac{\log\log n}{n}}}}\left|U_{n}(t)-U_{n}(t')-u(t)(t-t')\right|=o(n^{-\frac{1}{2}-\frac{\gamma}{8}}\log n).
\end{equation*}
\end{lem}

\begin{proof} As a consequence of Assumption \ref{ass3} and $\gamma<1$
\begin{equation*}
\sup_{\substack{t,t'\in[C_1,C_2]:\\ |t-t'|\leq C\sqrt{\frac{\log\log n}{n}}}}\left|U(t)-U(t')-u(t)(t-t')\right|=o(n^{-\frac{1}{2}-\frac{\gamma}{8}}\log n),
\end{equation*}
so it suffices to show that
\begin{equation*}
\sup_{\substack{t,t'\in[C_1,C_2]:\\ |t-t'|\leq C\sqrt{\frac{\log\log n}{n}}}}\left|U_{n}(t)-U_{n}(t')-(U(t)-U(t'))\right|=o(n^{-\frac{1}{2}-\frac{\gamma}{8}}\log n).
\end{equation*}
Without loss of generality, we can assume that $U(t)=t$, otherwise we use the same transformation as in the proof of Proposition \ref{pro1} and study the kernel function $h(x,y,U^{-1}(p))$. Note that in this case, we can consider the supremum over $[0,1]$. Furthermore, we will consider only the case $C=1$, we will prove
\begin{equation*}
K_n:=\sup_{\substack{t,t'\in[0,1]:\\ |t-t'|\leq \sqrt{\frac{\log\log n}{n}}}}\left|U_{n}(t)-U_{n}(t')-(t-t')\right|=o(n^{-\frac{1}{2}-\frac{\gamma}{8}}\log n).
\end{equation*}
For $l\in\N$, let $k=k_l=C2^{\lfloor\frac{1}{2}(l-\log\log l)\rfloor}$, so that for all $n=2^{l-1},\ldots,2^l-1$, we have that $\sqrt{\frac{\log\log n}{n}}\leq\frac{1}{k_l}\leq C\sqrt{\frac{\log\log n}{n}}$. We define for $r=0,\ldots,k_l$ the real numbers $t_{r,l}:=\frac{r}{k_l}$. Clearly 
\begin{multline*}
K_n\leq 2\max_{r=1,\ldots,k}\sup_{t,t'\in[t_{r-1,l},t_{r,l}]}\left|U_{n}(t)-U_{n}(t')-(t-t')\right|\\
\leq 4\max_{r=1,\ldots,k}\sup_{t\in[t_{r-1,l},t_{r,l}]}\left|U_{n}(t)-U_{n}(t_{r-1,l})-(t-t_{r-1,l})\right|.
\end{multline*}
Now chose $m=m_l\in\N$ such that $m_lk_l\approx2^{(\frac{1}{2}+\frac{\gamma}{8})l}$. So for all $n=2^{l-1},\ldots,2^l-1$ and some constants $C$, $C'$, we have that $Cn^{-\frac{1}{2}-\frac{\gamma}{8}}\leq\frac{1}{k_lm_l}\leq C'n^{-\frac{1}{2}-\frac{\gamma}{8}}$. We define for $r=1,\ldots,k_l$ and $r^\star=0,\ldots,m_l$ the real numbers $t^\star_{r^\star,r,l}=t_{r,l}+\frac{r^\star}{k_lm_l}$. As $U_n$ and $U$ are non-decreasing, we have for $t\in(t^\star_{r^\star-1,r,l},t^\star_{r^\star,r,l})$
\begin{multline*}
\left|U_n(t)-U_n(t_{r-1,l})-(t-t_{r-1,l})\right|\\
\leq\max\left\{\left|U_n(t^\star_{r^\star,r,l})-U_n(t_{r-1,l})-(t-t_{r-1,l})\right|\right.,\\
\left.\left|U_n(t^\star_{r^\star-1,r,l})-U_n(t_{r-1,l})-(t-t_{r-1,l})\right|\right\}\\
\leq\max\left\{\left|U_n(t^\star_{r^\star,r,l})-U_n(t_{r-1,l})-(t^\star_{r^\star,r,l}-t_{r-1,l})\right|\right.,\\
\left.\left|U_n(t^\star_{r^\star-1,r,l})-U_n(t_{r-1,l})-(t^\star_{r^\star-1,r,l}-t_{r-1,l})\right|\right\}+|t^\star_{r^\star,r,l}-t^\star_{r^\star-1,r,l}|,
\end{multline*}
and consequently
\begin{multline*} K_n\leq4\max_{r=1,\ldots,k}\max_{r^\star=1,\ldots,m}\left|U_{n}(t^\star_{r^\star,r,l})-U_{n}(t_{r-1,l})-(t^\star_{r^\star,r,l}-t_{r-1,l})\right|\\
+4\max_{r=1,\ldots,k}\max_{r^\star=1,\ldots,m}|t^\star_{r^\star,r,l}-t^\star_{r^\star-1,r,l}|\\
\leq8\max_{r=1,\ldots,k}\max_{r^\star=1,\ldots,m}\left|\frac{1}{n}\sum_{1\leq i\leq n} h_1(X_{i},t^\star_{r^\star,r,l})-\frac{1}{n}\sum_{1\leq i\leq n} h_1(X_{i},t_{r-1,l})\right|\\
+4\max_{r=1,\ldots,k}\max_{r^\star=1,\ldots,m}\left|\frac{2}{n(n-1)}\left(\sum_{1\leq i<j\leq n}h_2(X_{i
},X_{j},t^\star_{r^\star,r,l})-\sum_{1\leq i<j\leq n}h_2(X_{i},X_{j},t_{r-1,l})\right)\right|\\
+4\max_{r=1,\ldots,k}\max_{r^\star=1,\ldots,m}|t^\star_{r^\star,r,l}-t^\star_{r^\star-1,r,l}|.
\end{multline*}
By our construction of the numbers $t^\star_{r^\star,r,l}$, we have that $t^\star_{r^\star,r,l}-t^\star_{r^\star-1,r,l}=\frac{1}{k_lm_l}$ and obtain for all $n=2^{l-1},\ldots,2^l-1$
\begin{multline*}
\max_{r=1,\ldots,k}\max_{r^\star=1,\ldots,m}|t^\star_{r^\star,r,l}-t^\star_{r^\star-1,r,l}|\leq\sup_{t\in[C_1,C_2]}u(t)2^{-(\frac{1}{2}-\frac{\gamma}{4})l}\\
\leq Cn^{-\frac{1}{2}-\frac{\gamma}{8}} =o(n^{-\frac{1}{2}-\frac{\gamma}{8}}\log n).
\end{multline*}
With the help of Lemma \ref{lem3}, it follows that
\begin{multline*}
\max_{r=1,\ldots,k}\max_{r^\star=1,\ldots,m}\left|\frac{2}{n(n-1)}\left(\sum_{1\leq i<j\leq n}h_2(X_{i},X_{j},t^\star_{r^\star,r,l})-\sum_{1\leq i<j\leq n}h_2(X_{i
},X_{j},t_{r-1,l})\right)\right|\\
\leq \frac{4}{n(n-1)}\sup_{t\in\R}\left|\sum_{1\leq i<j\leq n}h_{2}\left(X_{i},X_{j},t\right)\right|=o\left(n^{-\frac{1}{2}-\frac{\gamma}{8}}\right).
\end{multline*}
Furthermore, we have for the linear part by Lemma \ref{lem1} and Corollary 1 of M\'oricz \cite{mori} (which gives moment bounds for the maximum other multidimensional partial sums)
\begin{multline*}
E\left[\left(\max_{n=2^{l-1},\ldots,2^{l}-1}\max_{r=1,\ldots,k}\max_{r^\star=1,\ldots,m}\left|\sum_{i=1}^n h_1(X_{i},t^\star_{r^\star-1,r,l})-\sum_{i=1}^n h_1(X_{i},t_{r-1,l})\right|\right)^4\right]\\
\leq\sum_{r=1}^kE\left[\left(\max_{n=2^{l-1},\ldots,2^{l}-1}\max_{m_1=1,\ldots,m}\left|\sum_{i=1}^n\sum_{r^\star=1}^{m_1}\left( h_1(X_{i},t^\star_{r^\star,r,l})- h_1(X_{i},t^\star_{r^\star-1,r,l})\right)\right|\right)^4\right]\\
\leq Ck2^{2l}l^2\left(\sqrt{\frac{\log l}{2^l}}\right)^{1+\gamma}=Cl^2(\log{l})^{\frac{\gamma}{2}}2^{(2-\frac{\gamma}{2})l},
\end{multline*}
as $E|h_1(X_{i},t)- h_1(X_{i},t')|\leq|t-t'|$ and by our construction $t^\star_{m,r,l}-t^\star_{0,r,l}=t_{r+1,l}-t_{r,l}=\frac{1}{k_l}\leq C\sqrt{\frac{\log l}{2^l}}$. So we can conclude that for any $\epsilon>0$
\begin{multline*}
\sum_{l=1}^\infty P\left[\max_{n=2^{l-1},\ldots,2^{l}-1}\max_{r\leq k}\max_{r^\star\leq m}\left|\sum_{i=1}^n \left(h_1(X_{i},t^\star_{r^\star-1,r,l})- h_1(X_{i},t_{r-1,l})\right)\right|\geq \epsilon 2^{\frac{1}{2}-\frac{\gamma}{8}l}l\right]\\
\leq C\sum_{l=1}^\infty\frac{2^{(2-\frac{\gamma}{2})l}l^2(\log{l})^{\frac{\gamma}{2}}}{\epsilon^4l^42^{(2-\frac{\gamma}{2})l}} =C\sum_{l=1}^\infty\frac{(\log l)^{\frac{\gamma}{2}}}{l^2}<\infty.
\end{multline*}
With the Borel Cantelli lemma, it follows that
\begin{equation*}
\max_{r=1,\ldots,k}\max_{r^\star=1,\ldots,m}\left|\sum_{1\leq i\leq n} h_1(X_{i},t^\star_{r^\star,r,l})-\sum_{1\leq i\leq n} h_1(X_{i},t_{r-1,l})\right|=o(n^{\frac{1}{2}-\frac{\gamma}{8}}\log n)
\end{equation*}
almost surely and finally
\begin{equation*}
\max_{r=1,\ldots,k}\max_{r^\star=1,\ldots,m}\left|\frac{1}{n}\sum_{1\leq i\leq n} h_1(X_{i},t^\star_{r^\star,r,l})-\frac{1}{n}\sum_{1\leq i\leq n} h_1(X_{i},t_{r-1,l})\right|=o(n^{-\frac{1}{2}-\frac{\gamma}{8}}\log n).
\end{equation*}
\end{proof}

\section{Proof of Main Results}

In all our proofs, $C$ denotes a constant and may have different values from line to line.

\begin{proof}[Proof of Theorem \ref{theo1}] We use the Hoeffding decomposition
\begin{equation*}
U_n\left(t\right)=U\left(t\right)+\frac{2}{n}\sum_{i=1}^{n}h_{1}\left(X_{i},t\right)+\frac{2}{n\left(n-1\right)}\sum_{1\leq i<j\leq n}h_{2}\left(X_{i},X_{j},t\right).
\end{equation*}
Let $K$ be a Gaussian process as in Proposition \ref{pro1}. Then
\begin{multline*}
\sup_{\substack{t\in\R\\s\in[0,1]}}\frac{1}{\sqrt{n}}\left|\lfloor ns\rfloor (U_{\lfloor ns\rfloor}(t)-U(t))-K(t,ns)\right|\\
\leq \sup_{\substack{t\in\R\\s\in[0,1]}}\frac{1}{\sqrt{n}s}\left|\left(2\sum_{1\leq i\leq ns}h_1(X_i,t)-K(t,ns)\right)\right|+\sup_{\substack{t\in\R\\s\in[0,1]}}\frac{1}{n^{\frac{3}{2}}s}\left|\sum_{1\leq i<j\leq ns}h_{2}\left(X_{i},X_{j},t\right)\right|\\
=O(\log^{-\frac{1}{3840}}n),
\end{multline*}
as by Lemma \ref{lem3}, we have
\begin{multline*}
\sup_{\substack{t\in\R\\s\in[0,1]}}\frac{1}{n^{\frac{3}{2}}s}\left|\sum_{1\leq i<j\leq ns}h_{2}\left(X_{i},X_{j},t\right)\right|\\
\leq n^{-\frac{\gamma}{8}}\sup_{\substack{t\in\R\\n'=1,\ldots n}}\frac{1}{(n')^{\frac{3}{2}-\frac{\gamma}{8}}}\left|\sum_{1\leq i<j\leq n'}h_{2}\left(X_{i},X_{j},t\right)\right|=O(n^{-\frac{\gamma}{8}}).
\end{multline*}
\end{proof}

\begin{proof}[Proof of Theorem \ref{theo2}] To simplify the notation, we will without loss of generality assume that $U(p)=p=t_p$ on the interval $I$. In the general case, one has to change the function $h(x,y,t)$ to $h(x,y,U^{-1}(t))$, as $Eh(X,Y,U^{-1}(p))=U(U^{-1}(p))=p$. The related empirical $U$-process $U_n\circ U^{-1}$, we have
\begin{multline*}
 R_n(p)=U_n^{-1}(p)-U^{-1}(p)-\frac{p-U_n(U^{-1}(p))}{u(t_p)}\\
 =\frac{1}{u(t_p)}\left((U_n\circ U^{-1})^{-1}(p)-p-(p-U_n\circ U^{-1}(p))\right)+o((U_n^{-1}(p)-U^{-1}(p))^{\frac{5}{4}}),
\end{multline*}
so Assumption \ref{ass3} guarantees that $R_n(p)$ is only blown up by a constant because of this transformation. If $U(p)=p=t_p$, then we can write $R_n(p)$ as
\begin{multline*}
 R_n(p)=U_n^{-1}(p)-t_p+U_n(t_p)-p\\
=\left(U_n^{-1}(p)-U_n^{-1}(U_n(t_p))+U_n(t_p)-p\right)+\left(U_n^{-1}(U_n(t_p))-t_p\right)
\end{multline*}
Applying Lemma \ref{lem5} and Lemma \ref{lem4} with $F=U_n$, $c=n^{-\frac{1}{2}-\frac{\gamma}{8}}\log n$ and $l=C\sqrt{\frac{\log\log n}{n}}$, we obtain
\begin{equation*}
\sup_{\substack{p,p'\in I:\\ |p-p'|\leq C\sqrt{\frac{\log\log n}{n}}}}\left|U_{n}^{-1}(p)-U_{n}^{-1}(p')-(p-p')\right|=o(n^{-\frac{1}{2}-\frac{\gamma}{8}}\log n).
\end{equation*}
almost surely. By Corollary \ref{cor1} we have that $\sup_{t\in[C_1,C_2]}\left(U_n(t_p)-p\right)\leq C\sqrt{\frac{\log\log n}{n}}$ almost surely, it follows that
\begin{multline*}
 \sup_{p\in I}\left|U_{n}^{-1}(p)-U_{n}^{-1}(U_{n}(t_p))+U_{n}(t_p)-p\right|\\
\leq\sup_{\substack{p,p'\in I:\\ |p-p'|\leq C\sqrt{\frac{\log\log n}{n}}}}\left|U_{n}^{-1}(p)-U_{n}^{-1}(p')-(p-p')\right|=o(n^{-\frac{1}{2}-\frac{\gamma}{8}}\log n)
\end{multline*}
almost surely. It remains to show the convergence of $U_n^{-1}(U_n(t_p))-t_p$. For every $\epsilon>0$ by the definition of the generalized inverse, $U_n^{-1}(U_n(t_p))-t_p>\epsilon n^{-\frac{1}{2}-\frac{\gamma}{8}}\log n$ only if $U_n(t_p+\epsilon n^{-\frac{1}{2}-\frac{\gamma}{8}}\log n)<U_n(t_p)$ and $U_n^{-1}(U_n(t_p))-t_p\leq-\epsilon n^{-\frac{1}{2}-\frac{\gamma}{8}}\log n$ only if $U_n(t_p-\epsilon n^{-\frac{1}{2}-\frac{\gamma}{8}}\log n)\geq U_n(t_p)$. So we can conclude that
\begin{multline*}
P\left[\sup_{p\in I}|U_{n}^{-1}(U_{n}(t_p))-t_p|> \epsilon n^{-\frac{1}{2}-\frac{\gamma}{8}}\log n\ \ \text{i.o.}\right]\\
\leq P\left[ \sup_{t\in[C_1,C_2-\epsilon n^{-\frac{1}{2}-\frac{\gamma}{8}}\log n]}U_{n}(t+\epsilon n^{-\frac{1}{2}-\frac{\gamma}{8}}\log n)-U_{n}(t)\leq0\ \ \text{i.o.}\right]\\
\leq P\left[ \sup_{\substack{t,t'\in[C_1,C_2]\\ |t-t'|=\epsilon n^{-\frac{1}{2}-\frac{\gamma}{8}}\log n}}\left|U_{n}(t)-U_{n}(t')-(U(t)-U(t'))\right|\geq|U(t)-U(t')|\ \ \text{i.o.}\right]\\
\leq P\left[ \sup_{\substack{t,t'\in[C_1,C_2]\\ |t-t'|\leq\epsilon n^{-\frac{1}{2}-\frac{\gamma}{8}}\log n}}\left|U_{n}(t)-U_{n}(t')-(U(t)-U(t'))\right|\geq\frac{\epsilon \log n}{n^{\frac{1}{2}+\frac{\gamma}{8}}\inf_{t\in[C_1,C_2]}u(t)}\ \ \text{i.o.}\right]\\
=0,
\end{multline*}
where the last line is a consequence of Lemma \ref{lem5}. We have proved that $\sup_{p\in I}|R_n(p)|=o(n^{-\frac{1}{2}-\frac{\gamma}{8}}\log n)$, and can finally conclude that
\begin{multline*}
\frac{n^{\frac{\gamma}{8}}}{\log n}\sup_{\substack{p\in I\\s\in[0,1]}}\frac{\lfloor ns\rfloor}{\sqrt{n}}|R_{\lfloor ns\rfloor}(p)|\\
\leq \sup_{n'\leq\sqrt{n}}(\frac{n'}{n})^{\frac{1}{2}-\frac{\gamma}{8}}\frac{\log n'}{\log n}\frac{n'^{\frac{1}{2}+\frac{\gamma}{8}}}{\log n'}\sup_{p\in I}|R_{n'}(p)|+\sup_{\sqrt{n}\leq n'\leq n}\frac{n'^{\frac{1}{2}+\frac{\gamma}{8}}}{\log n'}\sup_{p\in I}|R_{n'}(p)|\\
\leq C n^{-\frac{1}{4}+\frac{\gamma}{16}}\sup_{n'\in\N}\sup_{p\in I}|R_{n'}(p)|+\sup_{n'\geq \sqrt{n}}\frac{n'^{\frac{1}{2}+\frac{\gamma}{8}}}{\log n'}\sup_{p\in I}|R_{n'}(p)|\rightarrow 0.
\end{multline*}
\end{proof}

\begin{proof}[Proof of Theorem \ref{theo3}] Define $K'(p,s):=-\frac{1}{u(t_p)}K(t_p,s)$, there $K$ is the Gaussian process introduced in Theorem \ref{theo1}. $K'$ is then a Gaussian process with covariance function
\begin{equation*}
EK'(p,s)K'(p',s')=\min\left\{s,s'\right\}\frac{1}{u(t_p)u(t_{p'})}\Gamma(t_p,t_{p'})
\end{equation*}
and by Theorem \ref{theo1} and Theorem \ref{theo2}
\begin{multline*}
\sup_{\substack{p\in I\\s\in[0,1]}}\frac{1}{\sqrt{n}}\left|\lfloor ns\rfloor (U^{-1}_{\lfloor ns\rfloor}(p)-t_p)-K'(p,ns)\right|\\
\leq \sup_{\substack{p\in I\\s\in[0,1]}}\frac{1}{\sqrt{n}}\left|\lfloor ns\rfloor (U_{\lfloor ns\rfloor}^{-1}(p)-t_p-\frac{p-U_n(t_p)}{u(t_p)}\right|\\
+\sup_{\substack{p\in I\\s\in[0,1]}}\frac{1}{\sqrt{n}}\frac{1}{u(t_p)}\left|\lfloor ns\rfloor (U_{\lfloor ns\rfloor}(t_p)-p)-K(t_p,ns)\right|\\
\leq\sup_{\substack{p\in I\\s\in[0,1]}}\frac{\lfloor ns\rfloor}{\sqrt{n}}|R_{\lfloor ns\rfloor}(p)|+\frac{1}{\inf_{p\in I}u(t_p)}\sup_{\substack{p\in I\\s\in[0,1]}}\frac{1}{\sqrt{n}}\left|\lfloor ns\rfloor (U_{\lfloor ns\rfloor}(t_p)-p)-K(t_p,ns)\right|\\ 
=O(\log^{-\frac{1}{3840}}n)
\end{multline*}
almost surely.

\end{proof}

\begin{proof}[Proof of Theorem \ref{theo4}] If $\sigma^2>0$, set
\begin{equation*}
B(s)=\frac{1}{\sigma}T(K'(\cdot,s))=\int_I J(p)K'(p,s)dp+\sum_{j=1}^d b_jU_n(p_j).
\end{equation*}
In the case $\sigma^2=0$, $B$ may be an arbitrary Brownian motion. As $J$ is a bounded function, $T$ is a linear and Lipschitz continuous functional (with respect to the supremum norm), so
\begin{multline*}
\sup_{s\in[0,1]}\frac{1}{\sqrt{n}}\left|\lfloor ns\rfloor (T(U_{\lfloor ns\rfloor}^{-1})-T(U^{-1}))-\sigma B(ns)\right|\\
=\sup_{s\in[0,1]}\frac{1}{\sqrt{n}}\left|T\left(\lfloor ns\rfloor(U_{\lfloor ns\rfloor}^{-1}-U^{-1})-K'(\cdot,ns)\right)\right|\\
\leq C \sup_{\substack{p\in I\\s\in[0,1]}}\frac{1}{\sqrt{n}}\left|\lfloor ns\rfloor (U^{-1}_{\lfloor ns\rfloor}(p)-t_p)-K'(p,ns)\right|=O(\log^{-\frac{1}{3840}}n).
\end{multline*}
It remains to show that $B$ is a Brownian motion. Clearly, $EB(s)=0$ for every $s\geq0$. By the linearity of $T$, $B$ is a Gaussian process with stationary independent increments. Furthermore
\begin{multline*}
E[B^2(s)]=\frac{1}{\sigma^2}\int_{\tilde{C_1}}^{\tilde{C_2}}\int_{\tilde{C_1}}^{\tilde{C_2}}\frac{E[K(t_p,s)K(t_q,s)]}{u(t_p)u(t_q)}J(p)J(q)dpdq\\
+\frac{1}{\sigma^2}2\sum_{j=1}^db_j\int_{\tilde{C_1}}^{\tilde{C_2}}\frac{E[K(t_{p_j},s)K(t_q,s)]}{u(t_{p_j})u(t_p)}J(p)dp+\frac{1}{\sigma^2}\sum_{i,j=1}^db_ib_j\frac{E[K(t_{p_i},s)K(t_{p_j},s)]}{u(t_{p_i})u(t_{p_j})}\\
=s.
\end{multline*} 
\end{proof}

\section*{Acknowledgement}

We are very grateful for the careful reading and helpful comments of an anonymous referee which lead to a substantial improvement of the paper. The research was supported by the {\em Studienstiftung des deutschen Volkes} (German Academic Foundation) and the DFG Sonderforschungsbereich 823 (Collaborative Research Center) {\em Statistik nichtlinearer dynamischer Prozesse}.

\small{

}

\begin{thebibliography}{xxx}
	\bibitem{andr}{\scshape D.W.K. Andrews}, Non-strong mixing autoregressive processes, {\slshape J. Appl. Probab.} {\bfseries 21} (1984) 930-934.
	\bibitem{arc2}{\scshape M.A. Arcones}, The law of the iterated logarithm for $U$-processes, {\slshape J. Multivariate Anal.} {\bfseries 47} (1993) 139-151.
	\bibitem{arc6}{\scshape M.A. Arcones}, On the law of the iterated logarithm for Gaussian processes, {\slshape J. Theoret. Probab.} {\bfseries 8} (1995) 877-903.
	\bibitem{arco}{\scshape M.A. Arcones}, The Bahadur-Kiefer representation for $U$-quantiles, {\slshape Ann. Stat.} {\bfseries 24} (1996) 1400-1422.
	\bibitem{arc4}{\scshape M.A. Arcones, E. Gin\'e}, Limit Theorems for $U$-processes, {\slshape Ann. Prob.} {\bfseries 21} (1993) 1494-1542.
		\bibitem{arc5}{\scshape M.A. Arcones, E. Gin\'e}, On the law of the iterated logarithm for canonical $U$-statistics and processes, {\slshape Stochastic Process. Appl.} {\bfseries 58} (1995) 217-245.
	\bibitem{arc3}{\scshape M.A. Arcones, B. Yu}, Central limit theorem for empirical and $U$-processes of stationary mixing sequences, {\slshape J. Theoret. Probab.} {\bfseries 7} (1997) .47-53.
	\bibitem{aron}{\scshape N. Aronszajn}, Theory of reproducing kernels, {\slshape Trans. Amer. Math. Soc.} {\bfseries 3} (1950) 337-404.
	\bibitem{babu}{\scshape G.J. Babu, K. Singh}, On deviations between empirical and quantile processes for mixing random variables, {\slshape J. Multivariate Anal.}, {\bfseries 8} (1978) 532-549.
	\bibitem{baha}{\scshape R.R. Bahadur}, A note on quantiles in large samples, {\slshape Ann. Math. Stat.} {\bfseries 37} (1966) 577-580.
	\bibitem{berk}{\scshape I. Berkes, W. Philipp}, An almost sure invariance principle for the empirical distribution function of mixing random variables, {\slshape Probab. Theory Related Fields} {\bfseries 41} (1977) 115-137.
	\bibitem{boro}{\scshape S. Borovkova, R. Burton, H. Dehling}, Limit theorems for functionals of mixing processes with applications to $U$-statistics and dimension estimation, {\slshape Trans. Amer. Math. Soc.} {\bfseries 353} (2001) 4261--4318.
	\bibitem{brad}{\scshape R.C. Bradley}, {\slshape Introduction to strong mixing conditions}, volume 1-3, Kendrick Press, Heber City (2007).
	\bibitem{chou}{\scshape J. Choudhury, R.J. Serfling}, Generalized order statistics, Bahadur representations, and sequential nonparametric fixed-width confidence intervals, {\slshape J. Statist. Plann. Inference} {\bfseries 19} (1988) 269-282.
	\bibitem{csor}{\scshape M. C\"org\H{o}, P. R\'ev\'esz} Strong approximations of the quantile process, {\slshape Ann. Stat.} {\bfseries 4} (1978) 882-894.
	\bibitem{deh3}{\scshape H. Dehling, M. Denker, W. Philipp}, The almost sure invariance principle for the empirical process of $U$-statistic structure, {\slshape Annales de l'I.H.P.} {\bfseries 23} (1987) 121-134.
	\bibitem{dehl}{\scshape H. Dehling, M. Wendler}, Central limit theorem and the bootstrap for $U$-statistics of strongly mixing data, {\slshape J. Multivariate Anal.}, {\bfseries 101} (2010) 126-137.
	\bibitem{deh2}{\scshape H. Dehling, M. Wendler}, Law of the iterated logarithm for $U$-statistics of weakly dependent observations, in: Berkes, Bradley, Dehling, Peligrad, Tichy (Eds): {\slshape Dependence in Probability, Analysis and Number Theory}, Kendrick Press, Heber City (2010).
	\bibitem{denk}{\scshape M. Denker, G. Keller}, Rigorous statistical procedures for data from dynamical systems, {\slshape J. Stat. Phys.} {\bfseries 44} (1986) 67-93.
	\bibitem{geer}{\scshape J.C. Geertsema}, Sequential confidence intervals based on rank test, {\slshape Ann. Math. Stat.} {\bfseries 41} (1970) 1016-1026.
	\bibitem{hans}{\scshape B.E. Hansen}, GARCH(1,1) processes are near epoch dependent, {\slshape Econom. Lett.} {\bfseries 36} (1991) 181-186.
	\bibitem{hoef}{\scshape W. Hoeffding}, A class of statistics with asymptotically normal distribution, {\slshape Ann. Math. Stat.} {\bfseries 19} (1948) 293-325.
	\bibitem{hofb}{\scshape F. Hofbauer, G. Keller}, Ergodic properties of invariant measures for piecewise monotonic transformations, {\slshape Math. Z.} {\bfseries 180} (1982) 119-142.
	\bibitem{kief}{\scshape J. Kiefer}, Deviations between the sample quantile process and the sample df, in: M.L. Puri (Ed): {\slshape Nonparametric Techniques in Statistical Inference} (1970).
	\bibitem{kie2}{}{\scshape J. Kiefer}, Skorohod embedding of multivariate RV's, and the Sample DF, {\slshape Probab. Theory Related Fields} {\bfseries 24} (1972) 1-35.
	\bibitem{kuli}{\scshape R. Kulik}, Bahadur-Kiefer theory for sample quantiles of weakly dependent linear processes, {\slshape Bernoulli} {\bfseries 13} (2007) 1071-1090.
	\bibitem{lai}{\scshape T.L. Lai}, Reproducing kernel Hilbert spaces and the law of the iterated logarithm for Gaussian processes, {\slshape Probab. Theory Related Fields} {\bfseries 29} (1974) 7-19.
	\bibitem{mori}{\scshape F. M\'oricz}, A general moment inequality for the maximum of the rectangular partial sums of multiple series, {\slshape Acta Math. Hung.} {\bfseries 43} (1983) 337-346.
	\bibitem{mull}{\scshape D.W. M\"uller}, On Glivenko-Cantelli convergence, {\slshape Probab. Theory Related Fields} {\bfseries 16} (1970) 195-210.
	\bibitem{phil}{\scshape W. Philipp}, A functional law of the iterated logarithm for empirical functions of weakly dependent random variables, {\slshape Ann. Prob.} {\bfseries 5} (1977) 319-350.
	\bibitem{rous}{\scshape P.J. Rousseeuw, C. Croux}, Alternatives to the median absolute deviation, {\slshape J. Amer. Stat. Soc.} {\bfseries 88} (1993) 1273-1283.
	\bibitem{ser2}{\scshape R.J. Serfling}, The law of the iterated logarithm for $U$-statistics and related von Mises statistics, {\slshape Ann. Math. Statist.} {\bfseries 42} (1971) 1794.
	\bibitem{serf}{\scshape R.J. Serfling}, Generalized L-, M-, and R-statistics, {\slshape Ann. Prob.} {\bfseries 12} (1984) 76-86.
	\bibitem{wend}{\scshape M. Wendler}, Bahadur representation for $U$-quantiles of dependent data, {\slshape J. Multivariate Anal.}, {\bfseries 102} (2011) 1064-1079.
	\bibitem{wu}{\slshape W.B. Wu}, On the Bahadur representation of sample quantiles for dependent sequences, {\slshape Ann. Stat.} {\bfseries 33} (2005) 1934-1963.
	\bibitem{wu2}{\slshape W.B. Wu}, Strong invariance principles for dependent random variables, {\slshape Ann. Prop.} {\bfseries 35} (2007) 2294-2320.
	\bibitem{yosh}{\scshape K. Yoshihara}, Limiting behavior of $U$-statistics for stationary, absolutely regular processes, {\slshape Probab. Theory Related Fields} {\bfseries 35} (1976) 237-252.
\end{thebibliography}
\end{document}